# Blackwell's Demon – Postdiction and Prediction in Random Walks

By James D. Stein (California State University, Long Beach)

**Abstract** – Maxwell's Demon is a mythical being, first described by the physicist James Clerk Maxwell (although named Maxwell's Demon by Lord Kelvin). Maxwell used it in a thought experiment to potentially violate the Second Law of Thermodynamics by exploiting inhomogeneities existing in a statistically homogeneous system. Blackwell's Demon, making (as far as is known) its first appearance in this paper, illustrates a counterintuitive situation occurring in a random walk variation of the Two Envelope problem[1], that it is possible under restrictive conditions to predict with success probability > 1/2 the direction of a random walk generated by the flip of a fair coin. Like Maxwell's Demon, Blackwell's Demon operates by exploiting inhomogeneities that exist in a statistically homogeneous system. Maxwell's Demon achieves its results by knowing when a molecule is moving rapidly and when it is not. Blackwell's Demon achieves its results by knowing when a prediction strategy is successful and when it is not. At the time Maxwell proposed his Demon, it confronted a technological Everest – the ability to open and close a gate permitting the passage of a single molecule, and the ability to gauge the speed of an approaching molecule. Blackwell's Demon merely has to turn a light on and off, conduct visual observations and keep simple statistical records.

It should be noted that the analysis in this paper does not demonstrate the ability to predict the flip of a fair coin *ab initio* with success probability > 1/2, as it is necessary to embed the fair coin in an environment of some complexity in order to achieve this result.

**Introduction – Maxwell's Demon**

Maxwell imagined a canister of gas at thermal equilibrium, incapable of doing any work because there is no heat difference to convert to work. In the center of the canister is a gate, which can be opened and closed by the Demon, allowing the passage of a single molecule from one side to the other. Although the gas is at thermal equilibrium, some molecules are moving faster than others, as the temperature is an average obtained from the motion of all the molecules. When a hotter-than-average molecule on side A of the container approached the gate, the Demon opened it, allowing the molecule to pass to side B. Once the molecule had passed to side B, the Demon closed the gate. After some time, enough hot molecules were on side B to create a temperature difference and allow useful work to be done.

Maxwell's apparent intent in proposing the Demon was to demonstrate that the Second Law of Thermodynamics was a statistical regularity rather than an absolute certainty. Over the years, considerable debate ensued about the Demon and the potential violation of the Second Law. Physicists such as Leo Szilard and Rolf Landauer used the Demon to establish connections between information, entropy, reversible and irreversible thermodynamic processes.

This paper is not a paper about physics, but the background of Maxwell's Demon (and of course the proximity of the names Maxwell and Blackwell) prompted this paper's title. Additionally, as will be seen in the Remarks at the end of the article, there is considerable similarity between the two Demons, as indicated in the Abstract.

**Section I – Blackwell's Bet**

David Blackwell was one of the twentieth century's most influential statisticans, who also made significant contribution to probability theory and information theory. One of his ingenious decision strategies was transformed and renamed Blackwell's Bet by Len Wapner[3]. Wapner's version of the problem is described below.

There are two envelopes, a red envelope and a blue envelope, containing different sums of money. You select one at random, open it, and count the money. You are then offered the opportunity to either keep the money in that envelope, or open the other envelope and keep the money in that envelope. The goal, of course, is to end up with the larger amount of money.

Blackwell s contribution was to suggest that after opening the first envelope, the person then select a random number from any probability distribution whatsoever. Assume that the amount of money in the opened envelope is m, and the random number is r. Blackwell's strategy was to keep the money in the opened envelope if $r < m$, and the money in the unopened envelope if $r > m$.

The computation is easy. Suppose the two amounts of money in the two envelopes are S and L, with $S < L$. Initially, one has an equal probability of selecting each envelope. Let p be the probability that the random number selected is $> L$, and q the probability that the random number selected is $< S$.

If the envelope containing S is selected (which occurs with probability 1/2), the correct decision will be made with probability 1-q, for a combined probability of ½(1-q) of selecting the envelope containing S and making the correct decision.

If the envelope containing L is selected (which occurs with probability 1/2), the correct decision will be made with probability 1-p, for a combined probability of ½(1-p) of selecting the envelope containing L and making the correct decision.

The probability of making the correct decision is therefore ½(1-q) + ½(1-p) = ½ + ½ (1-(p+q)), or ½ + ½ the probability that the random number selected lies between S and L. As long as the probability of the random number lying between S and L is non-zero, which can be assured by many well-known distributions, the probability of a successful decision is > ½.

**Section II – A Postdiction Made with Help from Blackwell's Demon**

Maxwell's Demon has a job to do – opening and closing a gate based on information about an approaching molecule. Blackwell's Demon has a similar assignment, to turn lights on and off on a railroad track.based on information about the train, to observe the location of a random variable on the track, and to keep records.

The railroad track is circular, with N equally spaced stations. The track has a light midway between each adjacent pair of stations. The Demon, while riding on the train, has a laptop enabling him to turn on and off these lights. Constraints on N will be described later.

The train moves via a random walk generated by the flip of a fair coin. If the coin lands heads, the train moves 1 station clockwise, if the coin lands tails, the train moves 1 station counterclockwise.

We assume the limiting distribution has been reached, where the train has an equal probability of being at any station. The Demon is a passenger on the train. He does not know where he is, but hears the conductor announce, "Next stop, Willoughby." The Demon, a great fan of classic TV, is reminded of the episode "A Stop at Willoughby"[4] in the *Twilight Zone* television series from the 1960s.

At the moment the conductor announces the next stop, the Demon lights a light approximately diametrically opposite the Willoughby station. The Demon's laptop lights the light, but supplies no indication about location, either of the current station or Willoughby. Because the limiting distribution has been reached, the Demon is equally likely to be one stop clockwise, or one stop counterclockwise, from Willoughby.

Although the Demon does not know whether the coin landed heads or tails, the Demon can use a variant of Blackwell's decision strategy to guess with probability > ½ whether the coin landed heads or tails.

Fig. 1 shows the picture.

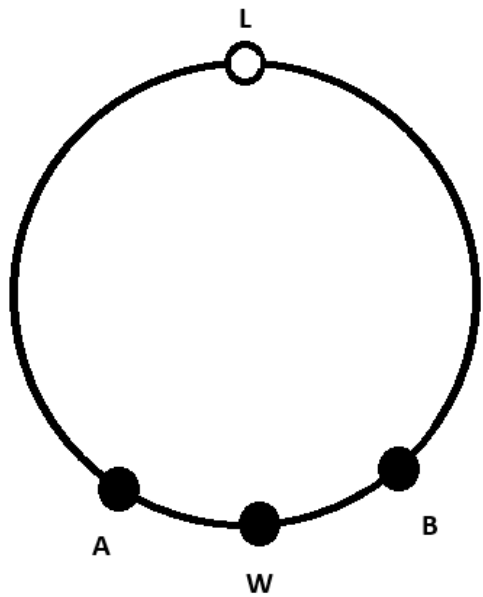

**Fig. 1**

Willoughby is represented by W, the light by L. Station A is one stop clockwise from W, station B is one stop counterclockwise from W. N is chosen sufficiently large to guarantee that L is on the opposite half of the circle as A, B and W (it is diametrically opposite in Fig. 1).

A variant of Blackwell's decision principle is implemented as follows. We will assume the track has circumference 1, and that a random location R on the track is chosen from the uniform distribution of locations on the track. Notice that whether W lies clockwise or counterclockwise from the train's location is equivalent to guessing whether the coin that determined the direction towards Willoughby landed heads or tails. From either A or B, there are two arcs from that location to L. R is on one of those two arcs. The Demon guesses that W lies on the arc in which R is encountered before L. If X is an arc on the circle, let len(X) denote the length of X.

The train is at A and B with equal probability. If the train is at A (probability = ½), the Demon will guess correctly if R is on the arc AWL, which occurs with probability len(AWL). The combined probability of being at A and guessing correctly is ½ len(AWL).

If the train is at B (probability = ½), the Demon will guess correctly if R is on the arc BWL, which occurs with probability len(BWL). The combined probability of being at A and guessing correctly is ½ len(BWL).

The combined probability of a correct guess is therefore ½ (len(AWL) + len(BWL)). As long as L is on the major arc joining A and B, as the choice of N ensures, the combined probability of a correct guess is ½ + ½ len(AB). Since len(AB) = 2/N, the probability of a correct guess is therefore ½ + 1/N.

This is a postdiction – the coin has been tossed, the destination decided. Even though the result of that 50-50 coin toss is unknown to the Demon, he can guess it correctly with probability ½ + 1/N.

The assumption that the light L is on the major arc joining A and B is needed to ensure that the probability of a successful postdiction is greater than ½. If the light is on the minor arc joining AB, the probability of a correct guess is easily seen to be ½ (len(AL) + len(BL)) = ½ len(AB) = 1/N.

If the light were in a fixed location that the Demon could not control, it would seem that the long-run average of the postdictions would be 1/2. Wherever the light is, there are only two destinations – the ones on either side of the light -- where the light is on the minor arc connecting the two potential stations of origin. With probability (N-2)/N, the destination would be such that the light is on the major arc joining the two stations on either side of the destination, and with probability 2/N the light is on the minor arc. The postdiction strategy described above in this case has a success probability (N-2)/N x (1/2 +1/N) + (2/N) x (1/N) = ½, as one might suspect.

Notice how the problem changes if Willoughby is the location of the train and not its destination. Then in Fig. 1, A and B are the two possible destinations, and a simple computation shows that the same strategy yields no advantage in predicting the toss of the coin which determines the next destination.

So, when the Demon can turn on a light guaranteeing the fidelity of Fig. 1, he can postdict the flip of the coin with probability > 1/2. But what if the coin that determines the next destination has not been flipped, so that the next destination has not yet been determined?

The Demon cannot turn on a light opposite the next destination, because that destination has not yet been determined. In the next section, we outline a strategy that makes it possible to guess the next destination – and thus the flip of the coin – with a probability greater than ½.

**Section III– Blackwell's Demon Makes a Prediction with Success Rate > 1/2**

At some time prior to steady-state, the Demon randomly turns on one of the lights (he could also turn on a specified light, it will make no difference to the following argument). That light remains on throughout the duration of the random walk. Its fixed location is not affected by the movement of the train.

Even though the coin has yet to be flipped, provided the train remains on the track, the next destination is $D_k$, which we abbreviate as D in the following diagram.

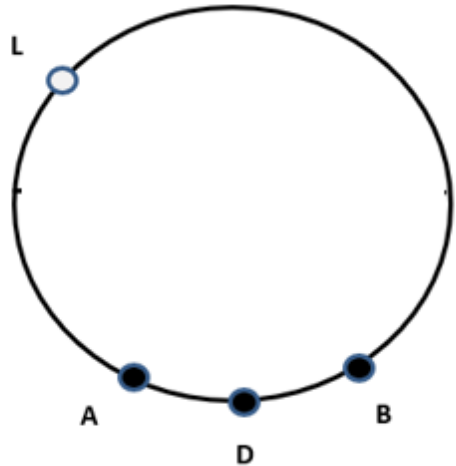

**Fig. 2**

The computations for this situation are identical to those for Fig. 1 when the Demon had no control over the location of the light. We will assume that the stations are numbered from 0 through N-1, and the light is located between stations 0 and 1. The overall success rate of the prediction strategy is ½ -- but when stations other than stations 0 or 1 are the next destination, the probability of a successful prediction is ½ + 1/N, and when stations 0 or 1 are the next destination, the probability of a successful prediction is 1/N. Define a "strong station" to be any station other than stations 0 or 1, and define a "weak station" to be either station 0 or station 1.

The Demon keeps track of his predictions, correlated by station number. Over time, the Demon will notice that at some stations – those with strong stations on either side – his prediction success rate is greater than ½, and at other stations his prediction success rate is less than ½. When this becomes statistically significant, the Demon modifies his predictions at those stations which have success rates less than ½, and simply guesses that the coin will land heads

(and so the train will proceed clockwise), This guess will be correct with probability ½. At all other stations, the probability of a correct guess is > ½, so the overall probability of a correct guess will be greater than ½.

The argument given here is easily adapted to the case of a random walk between a finite number of stations on a straight line. Thinking of this as a Markov chain with reflecting barriers, it is easily seen that if there are N stations on the line, the steady-state probabilities are 1/(2(N-1)) for the two end stations and 1/(N-1) for all other stations. The Blackwell pointer argument given in [5] using a uniformly-distributed random location yields a prediction success probability of (1/2 + 1/(N-1)) at the interior stations. Once again, a records-keeping Demon will guess correctly with probability > ½ at the interior stations and modify his guess at the other stations.

**Remarks**

The strategy outlined is not simply a matter of successfully predicting the toss of a fair coin more than half the time. There is a great deal of infrastructure in the current example, but it seems to provide an example of predicting a random walk with a greater success probability that one normally associates with it.

As with any counterintuitive conclusion, there is the possibility that an error has been made in the analysis. One possibility that might be argued is that because the destination has not yet been decided, one is not equally likely to be on either side of it. However, although the preceding analysis has been made for a single flip of the coin, it is true as an overall average – in the long run, whatever the destination is, the Demon will have been equally likely to be on either side of it, That is the only assumption on which the calculation for the probability of a successful guess is based. The average prediction success probability when a strong station is the next destination is ½ + 1/N, and when a weak station is the next destination, the average prediction success probability is 1/N.

It appears that, if the preceding analysis holds up, it is the result of the following combination of events.

(1) When this particular random walk achieves steady state, each station has an equal probability of being the next destination. Even if not all destinations are equally likely, the result can still hold. This is clear from an elementary continuity argument regarding the success probability as a function of the destination probabilities of the stations on either side. As long as there are stations at which the correct guess probability is greater than 1/2, statistical record-keeping should result in an ability to predict an outcome above what one would normally expect.

It is unclear whether this can occur for a simple random walk on the real line, as there is no steady state distribution. However, it may be possible to extend the results here to other random walks with steady state distributions.

(2) The presence of the light creates an inhomogeneity in the stations regarding the probability of a successful prediction.

(3) Statistical record-keeping can exploit the bias mentioned in (2).

**Similarities Between Maxwell's Demon and Blackwell's Demon**

**Existence of inhomogeneities in a statistically uniform system** – These inhomogeneities exist for Maxwell's Demon because the individual molecules move at different speeds. The inhomogeneities exist for Blackwell's Demon as a result of the existence of the light and the prediction strategy employed.

**Identification of inhomogeneities –** Maxwell's Demon is able to distinguish between molecules moving at different velocities.  Blackwell's Demon is able to distinguish between the differing prediction success rates at the various stations.

**Utilization of inhomogeneous components** – Maxwell's Demon uses a gate to separate the inhomogeneous components into two chambers with differing temperatures, allowing useful work to be performed by a heat engine.  Blackwell's Demon takes advantage of the knowledge of stations at which the adopted prediction strategy is performing "below par", and substitutes a better prediction strategy at these stations.

It does not appear that Maxwell believed that his Demon provided a practical way to violate the Second Law of Thermodynamics, but he would probably have been extremely interested to see the paths that resulted from his supposition.  Hopefully the discussion herein might also stimulate work under what situations an apparently counterintuitive successful prediction can arise.

**Acknowledgement**